\documentclass[journal]{ieeetran}
\usepackage{cite}
\usepackage{graphicx}
\usepackage{epstopdf}
\usepackage{amsmath}
\usepackage{cases}
\usepackage[caption=false,font=footnotesize]{subfig}
\usepackage{fixltx2e}
\usepackage{amssymb}
\usepackage{mathtools}
\usepackage{bm}
\usepackage{multirow}
\usepackage{color}
\usepackage[dvipsnames]{xcolor}
\usepackage{algorithm}
\usepackage{algorithmic}
\usepackage{amsthm}
\usepackage{mathrsfs}
\usepackage{textcomp,mathcomp}
\usepackage{booktabs}
\usepackage{threeparttable}

\newtheorem{remark}{Remark}
\usepackage{pifont}
\usepackage{threeparttable}

\newcommand{\tabincell}[2]{\renewcommand\arraystretch{0.9}\begin{tabular}{@{}#1@{}}#2\end{tabular}}

\begin{document}

\title{\LARGE{Harmonic Cancellation in Multi-Electrolyzer P2H Plants via Phasor-Modulated Production Scheduling}}

\author{
Yangjun~Zeng,~\IEEEmembership{Student Member,~IEEE},
Yiwei~Qiu,~\IEEEmembership{Member,~IEEE},
Li~Jiang,
Jie~Zhu,~\IEEEmembership{Student Member,~IEEE},
Yi~Zhou,~\IEEEmembership{Member,~IEEE},
Jiarong~Li,~\IEEEmembership{Member,~IEEE},
Shi~Chen,~\IEEEmembership{Member,~IEEE},
and
Buxiang~Zhou,~\IEEEmembership{Member,~IEEE}%

\vspace{-12pt}

\thanks{Financial support came from the National Key R\&D Program of China (2021YFB4000503) and the National Natural Science Foundation of China (52377116, 52577129, and 52307126). \emph{(Corresponding author: Yiwei Qiu)}}
\thanks{Y. Zeng, Y. Qiu, L. Jiang, J. Zhu, Y. Zhou, S. Chen, and B. Zhou are with the College of Electrical Engineering, Sichuan University, Chengdu 610065, China. \emph{(Corresponding author: Yiwei Qiu)}}
\thanks{J. Li is with the Harvard John A. Paulson School of Engineering and Applied Sciences, Harvard University, Cambridge 02138, USA.}%
}
\maketitle

\begin{abstract}
  Thyristor rectifiers (TRs) are cost-effective power supplies for hydrogen electrolyzers (ELZs) but introduce harmonic distortion that may violate grid codes. This letter proposes a self-governing harmonic mitigation strategy through coordinated operation of multiple ELZs in large power-to-hydrogen (P2H) plants. First, the harmonic model of TR-powered ELZs is derived, revealing a natural harmonic cancellation mechanism among them. Based on this, a system-level operation scheme based on phasor modulation is developed and integrated into plant scheduling. Case studies demonstrate that the proposed method reduces harmonic currents by 21.2\%--39.7\% and ensures grid-code compliance, with only a 0.25\% loss in hydrogen output, while increasing total revenue by over 21\% compared to production-oriented strategies.
\end{abstract}


\section{Introduction}
\label{sec:intro}

\IEEEPARstart{R}{enewable} power-to-hydrogen (ReP2H) offers promising pathways for green transition in the power and chemical sectors \cite{zeng2024scheduling}. As ReP2H projects expand in number and capacity, they employ multiple electrolyzers (ELZs), many powered by thyristor rectifiers (TRs) \cite{zeng2024scheduling} for their cost-effectiveness. However, the phase-controlled nature of TRs introduces harmonic distortion that may violate grid codes such as  IEEE~519 \cite{ieeestd} and GB/T~14549--93 \cite{GBT14549}; for example, a 12-pulse TR typically yields more than 6\% total harmonic distortion (THD) \cite{gao2024advanced}.

Existing harmonic mitigation approaches include multipulse rectification and active/passive power filters (APFs/PPFs), each with inherent tradeoffs. PPFs are inexpensive but unsuitable for fast-varying operating conditions \cite{das2004passive}. APFs and multipulse TRs are costly due to complex designs involving bidirectional PWM converters and phase-shifting transformers \cite{meng2021novel}. The industry requires low-cost harmonic mitigation solutions.

Several studies have addressed this challenge. Meng \textit{et al.} \cite{meng2021novel} proposed a hybrid rectifier, though its complexity hinders large-scale deployment. Yang \textit{et al.} \cite{yang2016enhanced} developed a phase-shifted current control scheme for motor drive systems, but it overlooks the on-off switching and load variations of ELZs.

In contrast to device-level approaches, this letter focuses on system-level solutions. It quantifies the coupling between harmonics and hydrogen production, identifies a tradeoff between harmonics and P2H efficiency, and proposes a harmonic mitigation method. The main contributions are:

1) Identification of a self-governing \textit{harmonic cancellation mechanism} in multi-ELZ systems, where adjusting electrolytic currents modulates harmonic phasors and reduces overall harmonic injections.

2) Formulation of a feasible region-based harmonic \textit{mitigation paradigm} and its integration into plant scheduling.

\begin{figure}[t]
  \centering
  \includegraphics[width=3.44in]{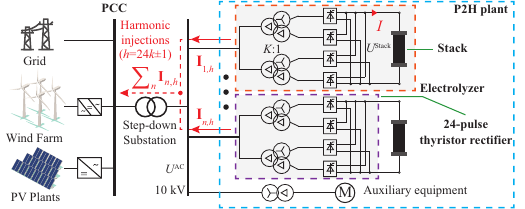}
  \caption{Schematic of the ReP2H system with TR-powered ELZs.}
  \label{fig:system}
\end{figure}

\section{Harmonic Model and Mitigation Strategy}

As shown in Fig. \ref{fig:system}, we assume each ELZ is independently powered by a 24-pulse rectifier\footnote{The proposed method also applies to TRs with different pulse numbers, as they can all have their phase modulated by the firing angle.}
without producing circulating currents \cite{zeng2024scheduling, li2024two},
which generates ($24k \pm 1$)th harmonic currents. The harmonic phasors from all ELZs are superposed at the point of common coupling (PCC) before being injected into the grid, which must satisfy grid-code limits.

\subsection{Harmonic Model of TR-Powered ELZs}
\label{sec:harmmodel}

\begin{figure}[t]
  \centering
  \includegraphics[width=3.45in]{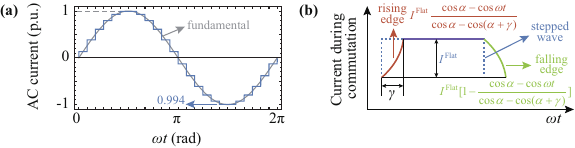}
  \caption{AC current waveform of the 24-TR. (a) The 24-step rectangular AC current. (b) The AC current during commutation.}
  \label{fig:24}
\end{figure}

Let the ELZs be rated at 5~MW, consisting of $N^\text{Cell}=350$ cells and operating within 2--7~kA.
The AC current has a 24-step rectangular waveform shown in Fig.~\ref{fig:24}(a). Due to line inductance, commutation takes time, as shown in Fig.~\ref{fig:24}(b) \cite{mohan2003power}. Considering the firing angle $\alpha$ and commutation overlap $\gamma$, a Fourier analysis yields the $h$th harmonic $\hat{\textbf{I}}_h(\alpha,\gamma) = {\textbf{I}_h}/{I_\text{1st}}$, where $\textbf{I}_h$ and $I_\text{1st}$ denote the $h$th harmonic and fundamental components, respectively. The analytical expression of $\hat{\textbf{I}}_h$ is not presented due to its complexity.

To associate  $\textbf{I}_h(\alpha,\gamma)$ with ELZ operation, the relation among $\alpha(I)$, $\gamma(I)$ and ${I}_\text{1st}(I)$ are established. Eqs. (\ref{eq:Ustack})--(\ref{eq:Ustack1}) relate the DC-side stack voltage $U^{\text{Stack}}$, commutation voltage drop $\Delta U$, and electrolytic current $I$, where $\Delta U$ is calculated by (\ref{eq:deltaU}), and ${I}_\text{1st}$ is determined via power conservation (\ref{eq:I1st}) \cite{li2024two}.
\begin{align}
	&  U^{\text{Stack}}=N^{\text{Cell}}[1.23 + 0.0776 I +0.07\log(1 + 54.86 I)],  \label{eq:Ustack}\\
	& U^{\text{Stack}}(I)=U^{\text{AC}}({2.4425}/{K})\cos\alpha-\Delta U, \label{eq:Ustack1}\\
	&\Delta U=U^{\text{AC}}\frac{2.4425}{K}\frac{\cos\alpha-\cos(\alpha+\gamma)}{2}=\frac{3}{\pi} X_\text{c} I, \label{eq:deltaU}\\
	& \sqrt{3}U^{\text{AC}}I_{\text{1st}} \big[{\cos\alpha+\cos(\alpha+\gamma)}\big]/{2}=U^{\text{Stack}}(I)I, \label{eq:I1st}
\end{align}
\noindent
where $U^{\text{AC}}=10$ kV is the AC-side bus voltage determined by the grid \cite{li2024two};
$K=16.67$ is the transformer turn ratio; and $X_\text{c}=0.0072~\Omega$ is the commutation reactance.
Solving \eqref{eq:Ustack}--\eqref{eq:I1st} simultaneously provides $\alpha(I)$, $\gamma(I)$, and $I_\text{1st}(I)$, which are substituted into $\hat{\textbf{I}}_h(\alpha,\gamma)$ to obtain $\textbf{I}_h(I)$.

\begin{figure}[t]
  \centering
  \includegraphics[width=3.45in]{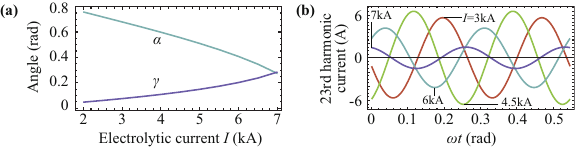}
  \caption{Variation of (a) $\alpha$ and $\gamma$, and (b) harmonic $\textbf{I}_\text{23rd}$ at various electrolytic current $I$.}
  \label{fig:harmonic23}
\end{figure}

Fig.~\ref{fig:harmonic23} shows $\alpha$, $\gamma$, and $\textbf{I}_\text{23rd}$ as functions of the stack current $I$. The phasors of the 23rd and 47th harmonics within the stack operation range [$\underline{I},\,\overline{I}$] are illustrated in Fig.~\ref{fig:gamma}. Higher-order harmonics are comparatively small and well below the grid-code thresholds \cite{GBT14549}, and thus not concerned here.

\subsection{Harmonic Cancellation Mechanism of Multiple ELZs}
\label{sec:harmcha}

\begin{figure}[t]
  \centering
  \includegraphics[width=3.35in]{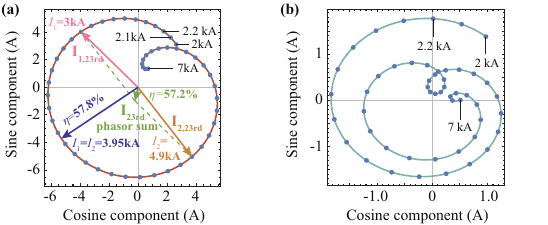}
  \caption{Harmonic current phasors as the electrolytic load current $I$ varies from 2 to 7 kA. (a) 23rd harmonic current. (b) 47th harmonic current.}
  \label{fig:gamma}
\end{figure}

\begin{figure}[t]
  \centering
  \includegraphics[width=3.45in]{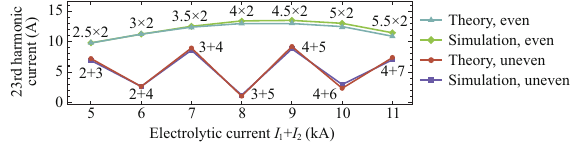}
  \caption{Validation of canceling $\textbf{I}_\text{23rd}$ with two ELZs via MATLAB/Simulink electromagnetic transient simulations under various combinations of $I$.}
  \label{fig:simulation}
\end{figure}

Fig.~\ref{fig:gamma}(a) shows that varying the electrolytic current $I$ modulates both the phase and magnitude of harmonic currents, enabling cancellation among multiple ELZs.
For instance, when two ELZs operate at $I_1=3$ kA and $I_2=4.9$ kA, their 23rd harmonic currents exhibit nearly equal magnitudes and opposite phases, resulting in a near-zero phasor sum.
Electromagnetic transient simulations in MATLAB/Simulink with two ELZs also validate the harmonic cancellation mechanism under different loading conditions, as shown in Fig. \ref{fig:simulation}.

However, this allocation is suboptimal for hydrogen production compared with the evenly distributed current of 3.95 kA (P2H efficiency $\eta$: 57.2\% $<$ 57.8\%; see the equimarginal principle in Appendix~A of our prior work \cite{zeng2024scheduling}).
Thus, a tradeoff arises between harmonic mitigation and energy efficiency, motivating the following system-level operational paradigm.

As the converter loss is approximately a quadratic function of $I$ \cite{li2024two, zheng2020discrete}, uneven current allocation may increase TR losses. However, the contribution of the quadratic term is small; therefore, the change in losses caused by distributing $I$ can be neglected. If higher accuracy is required, these losses can be incorporated into the power constraints of the scheduling model, as shown in Eq. (7) of \cite{zeng2024scheduling}.

\subsection{Operational Principle for Harmonic Mitigation}

This subsection formulates the harmonic feasible region of an $N$-ELZ P2H plant under GB/T 14549--93 \cite{GBT14549}, which requires the $h$th harmonic current at the PCC to remain within
\begin{align}
 \big| \sum\nolimits_{n=1}^N {\textbf{I}}_{n,h} \big|  \le \overline{I}_h=(S^{\text{sc}}/S^{\text{GB}})I_h^{\text{GB}},
\end{align}
\noindent where $S^{\text{sc}}$, $S^{\text{GB}}$ and $I_h^{\text{GB}}$ are the PCC short-circuit capacity, base short-circuit capacity, and corresponding harmonic limit. 

\begin{figure}[t]
  \centering
  \includegraphics[width=3.35in]{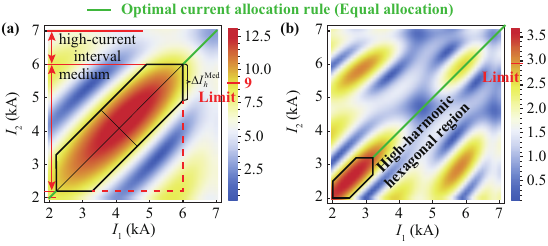}
  \caption{Amplitudes of the harmonic current of two ELZs. (a) 23rd harmonic current. (b) 47th harmonics current.}
  \label{fig:2ELZ}
\end{figure}

To satisfy these limits, the electrolytic loads of ELZs are coordinated to modulate their harmonic phasors. For clarity, two ELZs are first grouped to quantify harmonic limits as $2\overline{I}_h/N$ (with $N$ typically even in many projects \cite{zeng2025optimal}).
Numerical analysis of the 23rd and 47th harmonics shows that the infeasible region (IFR) is approximately a symmetric hexagon, as shown in Fig.~\ref{fig:2ELZ}(a).
Accordingly, the feasible region of the two-ELZ group can be expressed as
\begin{align}\label{eq:constraint23}
  \begin{cases}
  b_{n,h}^{\text{L}}+b_{n,h}^{\text{M}}+b_{n,h}^{\text{H}} =b_{n}^\text{On},~ h=23,25, \\
  I_{n}\leq\overline{I}_h^{\text{L}}b_{n,h}^{\text{L}}+\overline{I}_h^{\text{M}}b_{n,h}^{\text{M}}+\overline{I}_h^{\text{H}}b_{n,h}^{\text{H}}, \\
  I_{n}\geq\underline{I}_h^{\text{L}}b_{n,h}^{\text{L}}+\underline{I}_h^{\text{M}}b_{n,h}^{\text{M}}+\underline{I}_h^{\text{H}}b_{n,h}^{\text{H}}, \\
  \omega_h\leq b_{1,h}^{\text{M}},~\omega_h\leq b_{2,h}^{\text{M}},~\omega_h\geq b_{1,h}^{\text{M}}+b_{2,h}^{\text{M}}-1, \\
  \overline{N}_{h}-M(1-z_{h}) \leq \sum_n b_{n}^\text{On} < \overline{N}_{h}+M z_{h},  \\
  \left|I_{1}-I_{2}\right| \geq \Delta I_h^{\text{M}}-M(1-\omega_h)-M(1-z_h),
  \end{cases}
\end{align}
\noindent where
$b_{n}^\text{On}$ indicates whether the $n$th ELZ is active;
$b_{n,h}^{\text{L/M/H}}$ are binary variables indicating the ELZ current being in low/ medium/high-current intervals $\big[ \underline{I}_h^{\text{L/M/H}}, \overline{I}_h^{\text{L/M/H}} \big]$;
$\omega_h$ denotes whether both ELZs operate in the medium interval;
$\Delta I_h^{\text{M}}$ is the minimum electrolytic current difference ensuring harmonic compliance; $M$ is a large constant.

The shape of the IFR is explained as follows. With identical harmonic characteristics of the two ELZs, the IFR is inherently symmetric about the $45^\circ$ line, and harmonics decrease under uneven current allocation. As shown in Fig. \ref{fig:gamma}(a), the harmonic within 2--6 kA is roughly symmetric about 4 kA, and thus the IFR is roughly symmetric about $135^\circ$. To balance complexity and accuracy, the IFR is hence approximated by a hexagon.

\begin{remark}
  For multiple ELZ groups, the harmonic limit generalizes to $2\overline{I}_h/\sum_{n=1}^N b_n^\text{On}$, where idle groups enlarge the feasible region of active ones.
  When $\sum_{n=1}^N b_n^\text{On}\leq\overline{N}_h$, the $h$th harmonic remains compliant.
  Hence, the binary variable $z_h$ identifies whether harmonic mitigation is required.
\end{remark}

When the short-circuit capacity $S^{\text{sc}}$ changes, the harmonic limits and IFR change, but the shape of IFR remains a hexagon, allowing $S^{\text{sc}}$ variations to be handled by adjusting the parameters $\underline{I}_h^{\text{L/M/H}}, \overline{I}_h^{\text{L/M/H}}, \Delta I_h^{\text{M}}$, and $\overline{N}_h$ in (\ref{eq:constraint23}).

Although GB/T 14549 does not specify limits for the 47th and 49th harmonics, suppressing them remains beneficial. Restricting ELZs from operating in high-harmonic regions, as shown in Fig.~\ref{fig:2ELZ}(b), helps achieve this, and the same principle applies to current-distortion ratios defined in IEEE 519 \cite{ieeestd}.

\subsection{Integration into Plant Scheduling}

The proposed harmonic constraints (\ref{eq:constraint23}) can be easily incorporated into plant scheduling to jointly manage P2H efficiency and grid-code compliance.
A simplified scheduling model (\ref{eq:obj})--(\ref{eq:conslast}) is formulated here for clarity; in practice, these rules (\ref{eq:constraint23}) can be embedded in the comprehensive scheduling framework of our prior work \cite{zeng2024scheduling}, which also accounts for temperature dynamics, impurity limits of ELZs, and network power flow.
\begin{align}
	\max~ &\sum\nolimits_{t=1}^{T}\Big[\sum\nolimits_{n=1}^{N}(c^{\text{H}_2}Y_{n,t}^{\text{H}_2}-c^{\text{SU}}b_{n,t}^{\text{SU}})-c^{\text{G}}P_{t}^{\text{G}}\Big] \Delta t \label{eq:obj} \hspace{-0pt}\\
	\text{s.t.~~}&P_{t}^{\text{Rene}}+P_{t}^{\text{G}} \geq \textstyle\sum_{n=1}^{{N}} P_{n,t}^{\text{ELZ}},~P_{t}^{\text{G}}\geq0, \label{eq:power} \\
	& P_{n,t}^{\text{ELZ}}= P_{n,t}^{\text{Stack}}+(b_{n,t}^{\text{On}}+b_{n,t}^{\text{By}})P^{\text{Aux}}, \\
	&P_{n,t}^{\text{Stack}}=U^{\text{Stack}}(I_{n,t}) I_{n,t}, ~b_{n,t}^{\text{On}}\underline I\leq I_{n,t} \leq b_{n,t}^{\text{On}}\overline I,  \label{eq:PI} \\
	&Y_{n,t}^{\text{H}_2}=\eta^{\text{F}}N^{\text{Cell}}I_{n,t}/(2F), \label{eq:YH2}\\
	&b_{n,t}^{\text{On}}+b_{n,t}^{\text{By}}+b_{n,t}^{\text{Idle}}=1, \label{eq:logic}\\
	&b_{n,t}^{\text{On}}+b_{n,t}^{\text{By}}+b_{n,t-1}^{\text{Idle}}-1\leq b_{n,t}^{\text{SU}}, \label{eq:1}\\
	&b_{n,t-2}^{\text{Idle}}+b_{n,t}^{\text{Idle}}-b_{n,t-1}^{\text{Idle}} \geq 0, \label{eq:3}  \\
	&\text{ELZ current rules for harmonic mitigation (\ref{eq:constraint23})}, \label{eq:conslast}
\end{align}
\noindent
where the objective (\ref{eq:obj}) maximizes total revenue by balancing hydrogen output and operational costs;
(\ref{eq:power})--(\ref{eq:YH2}) describe power balance and hydrogen production; (\ref{eq:logic})--(\ref{eq:3}) represent the on-standby-idle transition of each ELZ;
$T$ and $\Delta t$ denote the scheduling horizon and step length; $c^{\text{H}_2}$, $c^{\text{G}}$, and $c^{\text{SU}}$ are the hydrogen, electricity, and startup cost; $P_{t}^{\text{Rene}}$ and $P_{t}^{\text{G}}$ are renewable and grid power; $P_{n,t}^{\text{ELZ}}$, $P_{n,t}^{\text{Stack}}$, and $P^{\text{Aux}}$ denote total, stack, and auxiliary power of each ELZ; $Y_{n,t}^{\text{H}_2}$ is the hydrogen output; $\eta^{\text{F}}$ and $F$ are the Faraday efficiency and constant; and $b_{n,t}^{\text{On}}$, $b_{n,t}^{\text{By}}$, $b_{n,t}^{\text{Idle}}$, and $b_{n,t}^{\text{SU}}$ indicate the operational states and startup action.

\section{Case Studies}
\label{sec:cases}

We compare the proposed method (\textbf{PM}) with two benchmarks:
\textbf{CM1} from \cite{qiu2023extended}, which optimizes hydrogen production without harmonic constraints; and
\textbf{CM2}, which follows CM1's on-off schedule but evenly distributes the load among ELZs, satisfying the harmonic constraint~(\ref{eq:constraint23}) by adjusting total load rather than using phasor modulation via load offsets.
A 2-ELZ system from \cite{li2024two} and a large-scale 20-ELZ system from \cite{zeng2024scheduling} are examined. 
Parameters are listed in Table~\ref{tab:para}.

\begin{table}[t]\scriptsize
  \renewcommand{\arraystretch}{0.95}
  \caption{Operational Parameters in the Case Study}
  \label{tab:para}
  \centering
  \begin{tabular}{c@{\hspace{4pt}}c@{\hspace{7pt}}c@{\hspace{4pt}}c@{\hspace{7pt}}c@{\hspace{4pt}}c}
  \toprule
  Parameter               & Value         & Parameter                & Value         & Parameter                & Value \\
   \midrule
   $c^{\text{H}_2}$     & 26 CNY/kg \cite{zeng2024scheduling}    &  $c^{\text{G}}$       &  0.6 CNY/kWh  &$c^{\text{SU}}$     & 1000 CNY \cite{zeng2024scheduling} \\
   $P^{\text{Aux}}$     &  0.5 MW    &$\underline I$, $\overline I$       & 2 kA, 7 kA         &  $T$,$\Delta t$     &  24, 1 h        \\
  \bottomrule
  \end{tabular}
\end{table}

\begin{figure}[t]
	\centering
	\includegraphics[width=3.25in]{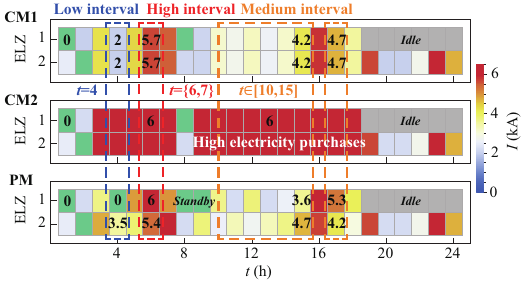}
	\caption{State transitions and load allocation in the 2-ELZ system under the conventional methods (CM1 and CM2) and the proposed method (PM).}
	\label{fig:rule}
\end{figure}

\begin{figure}[t]
	\centering
	\includegraphics[width=3.5in]{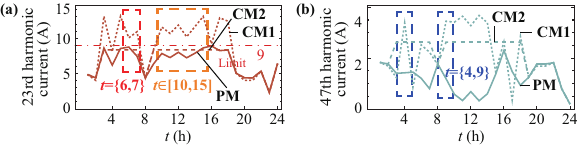}
	\caption{Comparison of (a) 23rd and (b) 47th harmonic currents in the 2-ELZ system under the different scheduling methods.}
	\label{fig:harmonicI}
\end{figure}

\textit{1) 2-ELZ Illustrative Case:}
In the 2-ELZ system, the rated voltage is 10~kV and the short-circuit capacity is $S^{\text{sc}}=200$~MVA.
Figs.~\ref{fig:rule} and~\ref{fig:harmonicI} compare the ELZ state transitions, load allocation, and harmonic currents under CM1, CM2, and PM.

During high-current intervals, $t=[6,7]$~h, PM mitigates minor harmonic violations by slightly uneven load allocation among ELZs. In medium-load periods, $t \in [10,15]$~h, CM1 exceeds the 23rd and 25th harmonic limits, whereas PM introduces small inter-ELZ current offsets to modulate phase differences, avoiding harmonic hotspots; see Fig.~\ref{fig:2ELZ}(a).
Simultaneously, the 47th and 49th harmonics are also suppressed.
At low loads, $t=4$ or $9$~h, rule~(\ref{eq:constraint23}) prevents both ELZs from operating at low current; thus, PM temporarily switches one unit to standby, effectively reducing harmonic injections.

Table~\ref{tab:comparison} summarizes the performance metrics.
CM1 yields the highest revenue but violates grid codes.
PM reduces average harmonic currents by 19.4\%--42.8\% relative to CM1, with only a 0.8\% revenue loss and 0.02\% hydrogen output loss.
In CM2, current is evenly allocated, as indicated by the green diagonal in Fig. 6(a), but harmonic limits force currents above 6 kA or below 2.2 kA, resulting in increased power purchase or renewable energy curtailment. Because grid electricity costs 600 CNY per MWh while producing hydrogen worth only 494 CNY, higher electricity purchases reduce the total revenue.
Compared with CM2, PM achieves 28\% higher revenue under similar harmonic levels, highlighting the effect of phasor modulation.

\begin{table}[t]\scriptsize
	\renewcommand{\arraystretch}{1.05}
	\caption{Performance Comparison in the 2-ELZ Case}
	\label{tab:comparison}
	\centering
	\begin{tabular}{c@{\hspace{4pt}}c@{\hspace{4pt}}c@{\hspace{4pt}}c}
		\toprule
		Method                                                    & \textbf{CM1}    & \textbf{CM2} & \textbf{PM}    \\
		\midrule
		\tabincell{c}{Revenue (10$^3$ CNY)}                       &   47.63        &  36.90      &  47.25          \\
		\tabincell{c}{Hydrogen output (kg)}                        &   1867.2       &  2734.27    & 1866.8         \\
		\tabincell{c}{Grid electricity purchase (MWh)}                &   1.53         &  56.99      &  2.14          \\
		\tabincell{c}{Grid codes compliance}                         & \ding{55}       & \ding{51}  & \ding{51}    \\
		\tabincell{c}{Average 23/25/47/49th \\harmonic currents (A)}   & 8.7/6.7/2.1/1.9         &7.0/\text{5.1}/2.2/2.0        &\text{6.7}/5.4/\text{1.2}/\text{1.1}     \\
		\bottomrule
	\end{tabular}
\end{table}

\textit{2) 20-ELZ Industrial Case:}
In the 20-ELZ system, the rated voltage is 220~kV and $S^{\text{sc}}=3,900$~MVA.
Table~\ref{tab:comparison1} compares CM1, CM2, and PM results.
Relative to CM1, PM achieves 21.2\%--39.7\% reductions in harmonic currents with only 0.25\% revenue and hydrogen output losses, while outperforming CM2 by over 21\% in revenue under comparable harmonic levels.
With higher system flexibility, large-scale plants benefit more significantly from coordinated harmonic management, validating the scalability and practical value of the proposed approach.

\begin{table}[t]\scriptsize
	\renewcommand{\arraystretch}{1.05}
	\caption{Performance Comparison in the 20-ELZ Case}
	\label{tab:comparison1}
	\centering
	\begin{tabular}{c@{\hspace{4pt}}c@{\hspace{4pt}}c@{\hspace{4pt}}c}
		\toprule
		Method                                                    & \textbf{CM1}    & \textbf{CM2}      & \textbf{PM}    \\
		\midrule
		\tabincell{c}{Revenue (10$^3$ CNY)}                       & 479.2         &  394.1          &  478.0   \\
		\tabincell{c}{Hydrogen output (kg)}                        & 18463         & 25200           &  18417  \\
		\tabincell{c}{Grid electricity purchase (MWh)}                &  1.37         &  435.23         &  1.35          \\
		\tabincell{c}{Grid codes compliance}                         & \ding{55}       & \ding{51}  & \ding{51}    \\
		\tabincell{c}{Average 23/25/47/49th \\harmonic current (A)}   & 3.3/2.8/0.7/0.7         &\text{2.6}/\text{1.8}/0.8/0.7        &\text{2.6}/2.1/\text{0.5}/\text{0.4}     \\
		\bottomrule
	\end{tabular}
\end{table}

\section{Conclusions}
\label{sec:conclusion}

This letter addresses the harmonic compliance problem in P2H plants by establishing a harmonic model that quantifies the coupling between electrolytic currents and harmonic distortion.
Building on this model, a system-level phasor modulation strategy is developed to coordinate multiple ELZs for self-governing harmonic cancellation. Rules are then formulated to integrate this mechanism into plant operation.

Case studies show that the dominant 23rd/25th harmonics arise in the medium-load range and the 47th/49th harmonics in the low-load range.
The proposed method effectively suppresses these harmonics, ensuring grid-code compliance with only 0.25\% reductions in hydrogen output and revenue. 

Future research will investigate online harmonic monitoring and optimal allocation of PPF/APFs and other resources to enable more grid-friendly large-scale P2H deployment.

\vfill
\break


\begin{thebibliography}{10}
\providecommand{\url}[1]{#1}
\csname url@samestyle\endcsname
\providecommand{\newblock}{\relax}
\providecommand{\bibinfo}[2]{#2}
\providecommand{\BIBentrySTDinterwordspacing}{\spaceskip=0pt\relax}
\providecommand{\BIBentryALTinterwordstretchfactor}{4}
\providecommand{\BIBentryALTinterwordspacing}{\spaceskip=\fontdimen2\font plus
\BIBentryALTinterwordstretchfactor\fontdimen3\font minus
  \fontdimen4\font\relax}
\providecommand{\BIBforeignlanguage}[2]{{%
\expandafter\ifx\csname l@#1\endcsname\relax
\typeout{** WARNING: IEEEtran.bst: No hyphenation pattern has been}%
\typeout{** loaded for the language `#1'. Using the pattern for}%
\typeout{** the default language instead.}%
\else
\language=\csname l@#1\endcsname
\fi
#2}}
\providecommand{\BIBdecl}{\relax}
\BIBdecl

\bibitem{zeng2024scheduling}
Y.~Zeng, Y.~Qiu, J.~Zhu, S.~Chen, B.~Zhou, J.~Li, B.~Yang, and J.~Lin,
  ``Scheduling multiple industrial electrolyzers in renewable {P2H} systems: A
  coordinated active-reactive power management method,'' \emph{IEEE Trans.
  Sustain. Energy}, vol.~16, no.~1, pp. 201--215, Jan. 2025.

\bibitem{ieeestd}
``{IEEE} recommended practice and requirements for harmonic control in electric
  power systems,'' \emph{IEEE Std 519-2014 (Revision of IEEE Std 519-1992)},
  pp. 1--29, 2014.

\bibitem{GBT14549}
``Quality of electric energy supply--{H}armonics in public supply network,''
  Standard GB/T 14549-93, National Standard of the People's Republic of China,
  1993, in Chinese.

\bibitem{gao2024advanced}
Y.~Gao, X.~Wang, and X.~Meng, ``Advanced rectifier technologies for
  electrolysis-based hydrogen production: A comparative study and real-world
  applications,'' \emph{Energies}, vol.~18, no.~1, p.~48, 2024.

\bibitem{das2004passive}
J.~Das, ``Passive filters-potentialities and limitations,'' \emph{IEEE Trans.
  Ind. Appl.}, vol.~40, no.~1, pp. 232--241, Jan. 2004.

\bibitem{meng2021novel}
X.~Meng, M.~Chen, M.~He, X.~Wang, and J.~Liu, ``A novel high power hybrid
  rectifier with low cost and high grid current quality for improved efficiency
  of electrolytic hydrogen production,'' \emph{IEEE Trans. Power Electron.},
  vol.~37, no.~4, pp. 3763--3768, Apr. 2022.

\bibitem{yang2016enhanced}
Y.~Yang, et al., ``Enhanced phase-shifted current
  control for harmonic cancellation in three-phase multiple adjustable speed
  drive systems,'' \emph{IEEE Trans. Power Deliv.}, vol.~32, no.~2, pp.
  996--1004, Apr. 2017.

\bibitem{li2024two}
J.~Li, et al., ``Two-layer
  energy management strategy for grid-integrated multi-stack power-to-hydrogen
  station,'' \emph{Appl. Energy}, vol. 367, p. 123413, Aug. 2024.

\bibitem{mohan2003power}
N.~Mohan, T.~M. Undeland, and W.~P. Robbins, \emph{Power Electronics:
  Converters, Applications, and Design}.\hskip 1em plus 0.5em minus 0.4em\relax
  John Wiley \& Sons, 2003.

\bibitem{zheng2020discrete}
J.~Zheng, et al., ``A discrete state event
  driven simulation based losses analysis for multi-terminal megawatt power
  electronic transformer,'' \emph{CES Trans. Elect. Machines Syst.}, vol.~4,
  no.~4, pp. 275--284, Dec. 2020.

\bibitem{zeng2025optimal}
Y.~Zeng, et al., ``Optimal
  investment portfolio of thyristor- and {IGBT}-based electrolysis rectifiers
  in utility-scale renewable {P2H} systems,'' \emph{IEEE Trans. Sustain.
  Energy}, 2025, early access.

\bibitem{qiu2023extended}
Y.~Qiu, et al.,
  ``Extended load flexibility of utility-scale P2H plants: Optimal production
  scheduling considering dynamic thermal and HTO impurity effects,''
  \emph{Renewable Energy}, vol. 217, p. 119198, 2023.

\end{thebibliography}
\end{document}